# 事故性溢油情景下船舱油液泄漏入水模拟方法

**摘要**：船舶溢油事故已成为海洋环境污染与生态损害的主要原因之一，运用模拟预测方法研究油液泄漏入水过程对油污清理以及应急救援工作有重大意义。对溢油模拟预测做了大量深入的研究，研究分析了溢油油液泄漏入水行为过程研究的发展过程以及面临的难题，对简化模拟模型、小孔射流模型和计算流体动力学模型这3类油液泄漏入水过程的优缺点进行了分析比较。在此基础上，指出了研究中存在的问题以及今后研究的重点，提出建立复杂船体模型、复杂海况条件相互叠加以及模拟计算方法耦合的新思路。

**关键字**：溢油污染；入水过程；复杂船体结构；环境因素

## Simulation Method Of Cabin Oil Leakage Into Water In Case Of Accidental Oil Spill

**Abstract:** Ship oil spill accident has become one of the main causes of Marine environmental pollution and ecological damage. It is of great significance to study the process of oil spill into water by using simulation prediction method for oil spill clean-up and emergency rescue work. A great deal of in-depth research on oil spill simulation and prediction is done. The development process of the research on the behavior process of oil spill into water is analyzed and the problems are put forward. The advantages and disadvantages of simplified simulation model, keyhole jet model and computational fluid dynamics model are analyzed and compared. On this basis, the problems existing in the research and the emphases of future research are pointed out, and a new idea of building complex hull model, superposition of complex sea conditions and coupling of simulation calculation methods is proposed.

## 引言

据国际油轮船东污染组织（ITOPF）的统计，1970-2020年全球范围内发生7吨及以上船舶溢油事故1845起，总溢油量约586万吨，低于7吨的溢油事故更是不计其数[1]。海上溢油事故频发不仅使自然环境、生态系统受到严重损害，经济蒙受损失，而且还严重危害人体健康，其造成的污染更会持续数十年[2]。为了对船舶溢油进行防控，世界上很多国家在防污法律法规的制定，溢油防污技术的开发，以及溢油模拟预测等方面做了大量的研究工作。其中，溢油模拟预测能够预测油类污染物在海上的行为过程，为油污清理以及生态环境损害评估提供科学支撑信息[3]。溢油模拟预测包括溢油油液泄漏入水过程模拟、油膜漂移扩散模拟、溢油风化过程模拟、溢油对生态损害模拟等等。油膜漂移扩散、溢油风化等方面的模拟做了大量研究，也取得了很多研究成果，但是在溢油油液泄漏入水过程模拟的研究较少，进展也十分缓慢，究其原因主要是该过程的准确模拟涉及到复杂船体结构建模和事故发生时的复杂海况条件耦合运算，模拟起来非常困难，目前业界普遍做法是对整个泄漏过程情景参数做不同程度的简化，如忽略了双层壳、货舱型式和泄漏口特征等，但这势必影响泄漏过程模拟结果

**作者简介**：顾乾乾（1996-），山东枣庄人，硕士，研究方向船舶与海洋工程。E-mail:1476480719@qq.com
**通信作者**：庄学强（1974-），福建泉州人，教授，博士，研究方向船舶安全监测与污染防治技术，13015917291@163.com

的准确性。船舶溢油事故中油液泄漏入水过程模拟结果是后续的溢油漂移、扩散、风化行为归宿模拟计算和溢油环境损害预测计算的重要输入源项，源项的准确性将会直接影响了后续的模拟预测的准确性，进而影响了应急措施的科学性和损害索赔的客观性。因此溢油油液泄漏入水过程模拟研究对整个溢油模拟以至于溢油应急技术的研究有十分重要的意义。

# 1 基于溢油量估算的泄漏过程极简化模拟

溢油量是衡量船舶溢油事故严重程度最重要的参数之一，我国《中华人民共和国水污染防治法》和《防治船舶污染海洋环境管理条例》等法律法规也专门规定了事故发生后溢油量与污染事故等级划分和行政处罚额度的对应关系。于是很长一段时间以来，人们经常把溢油泄漏行为描述简化处理为溢油量估算，这是一种只考虑泄漏结果，不考虑过程的极简化模拟方法。目前关于泄漏量估算方法有HFO（Heavy Fuel Oil；重燃油）变换量估算法、油膜监测估算法、光学估计法等等。

HFO变换量估算法是比较基础的估算方法，依据船舶日常油量记录数据，HFO总库存量变化与HFO设备运转消耗量的减少来推算出溢油量，公式为[4]：

$$G = Z - C - R \quad (1)$$

式中，$G$为船舶破损后估计的溢油量；$Z$为破损前船舶总储存油量；$C$为船舶记录总油量时刻至船舶破损时这期间船舶的总耗油量；$R$为船舶破损后油舱剩余油量。

溢油油膜监测估算法是在船舶污染事故发生后，通过现场观测、卫星遥感等多种监测技术获取油膜面积，根据波恩协议中油膜颜色与油膜厚度的对应关系获得油膜厚度，其海面溢油量计算公式为[5]：

$$G = \sum_{i=1}^{n} S_i \times H_i \times \rho \quad (2)$$

式中，$G$为溢油估计量；$S_i$为第$i$种颜色的油膜面积；$H_i$为第$i$种颜色油膜的厚度；$\rho$为油品密度；$n$为油膜颜色分区数。

光学估计法是利用粒子图像测速技术，通过对2个连续的视频进行观察，分析视频帧之间的距离得出溢油的瞬时速度或是平均速度再与破损面积、油品密度和溢油时间相乘便得到溢油估计量[6]，公式为：

$$G = S \times V \times \rho \times t \quad (3)$$

式中，$G$为溢油估计量；$S$为破损处横截面积；$V$为溢油流速；$\rho$为油品密度，$t$为溢油泄漏时间。

通过先进的海洋观测设备采集数据，简化处理后进行溢油量的估计，在实验研究中有很大的可行性，但是溢油事故常常发生在恶劣的海况条件下以及存在突发性，这使得数据采集存在极大的困难，在实际应用中很难适用。简化处理溢油速度也是目前常用的OILMAP等溢油软件前处理采用的方法，首先对溢油量进行估计，再推算出溢油持续的时间，将溢油过程简化为瞬时泄漏或者连续定常泄漏行为来进行溢油的后续预测。

HFO变换量估算法、油膜监测估算法仅仅关注了泄漏结果，只得到一个静态的泄漏量数据，而光学估计法、OILMAP等溢油软件只是更进一步考虑了时间因素，简单的将溢油过程简化为连续定常泄漏或者瞬时泄漏过程，都没有考虑泄漏过程中液位、溢油速度等参数的动态变化对溢油泄漏的影响。

# 2 基于伯努利方程的小孔射流模拟方法

DNV（Det Norske Veritas；挪威船级社）、ABS（American Bureau of Shipping；美国船级社）等著名机构建议采用孔洞射流模式来分析溢油流动状态以及溢油泄漏行为，孔洞射流模式是随液面高度和泄漏时间变化的泄漏源强模型，更加接近于船舶溢油实际泄漏行为[7]。船舶发生破损，穿孔处与外部大气或是海水相通造成液体流出，研究分析破孔位置在水线上和水线下造成的不同泄漏情况，如图1、图2所示。

破孔位置在水线上的泄漏行为较为简单，主要影响因素是受到液体重力，图1为油船侧壁破孔在水线之上横截面示意图。破孔面积用$A_h$表示，油舱的上表面积为$A_t$，破裂的流体速度的时变静水压头$h$是从破裂中心线到油舱油的自由表面的垂直距离，用时间标度$t_d$来表征油舱破裂后溢油的持续时间。对水线之上的泄漏进行理想化处理，这里假设了破孔面积相较于油舱较小，因此溢油流出的速度假设为定常流动，通过破孔的流出速度大小为$\sqrt{2gh_0}$，因此流出体积流率$\sqrt{2gh_0}A_h$乘以泄漏时间$t_d$等于泄漏的液体体积流量（其中$h_0$是$h$的初始值），得出：

$$Q = \sqrt{2gh_0} A_h t_d \quad (4)$$

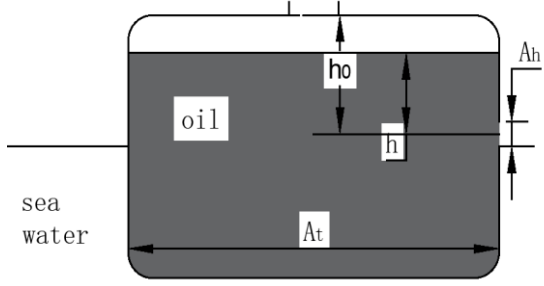

图 1 油舱侧壁水线上破孔示意图

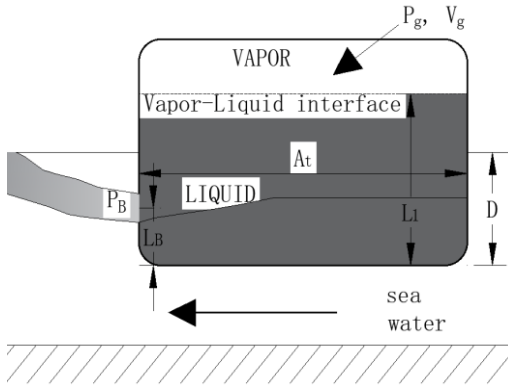

图 2 油舱水线下破孔示意图

破孔位置在水线下的泄漏行为较为复杂，涉及到的影响因素较多，采用孔洞射流模型，该模型由伯努力方程（Bernoulli Equation）推导而得，简单表述如下：

$$Q = C_D A_B \rho_L \left[ \frac{2(P_T - P_{out})}{\rho_L} + 2gH \right]^{\frac{1}{2}} \quad (5)$$

式中：$C_D$ 为泄漏口的流量系数；$A_B$ 为泄漏口的流通面积，$m^2$；$P_T$ 为油舱上方的压力，$Pa$；$P_{out}$ 破损口舷外压力（根据破损口位置，可能是大气压或者是水压），$Pa$；$\rho_L$ 为油液密度，$kg/m^3$；$H$ 为液面离泄漏孔中心距离，$m$；$Q$ 为油液泄漏率，$kg/s$。

式（5）是一个理想化的方程，仅仅考虑了泄漏过程中随着油液下降造成压力变化从而影响油液泄漏率，这与泄漏发生时的复杂海况条件不相符合。研究者在此基础上进一步研究了水线下静水压力以及船体运动造成的泄漏，开发了一种模型来确定由于船体运动或是水面波动引起的重力和压力变化而引起的泄漏，为两个排放阶段建立了分析模型：（1）初始阶段，持续到通过断裂处建立静水压力平衡，（2）由于船体运动改变静水压力平衡而产生的泄漏[8-9]。

第一阶段破孔处瞬时速度公式：

$$u_B(t) = -\left[ g\frac{A_B}{A_t} + \frac{P_g(0)V_g(0)}{(V_g)^2}\frac{A_B}{\rho_l} \right] C_D^2 t + u_B(0) \quad (6)$$

$$u_B(0) = C_D \sqrt{2(L_l - L_B)g - 2\frac{\rho_w}{\rho_l}(D - L_B)g} \quad (7)$$

第二阶段破孔处的瞬时速度公式：

$$u_B = C_D \sqrt{2ag(\frac{\rho_w}{\rho_l})\sin\frac{2\pi t_2}{T}} \quad (8)$$

式中：$L_l$ 为舱内液位；$L_B$ 为从船底至船体破孔处距离；$P_g$ 为舱内气体空间压力；$P_a$ 为大气压力；$P_B$ 为破孔处压力；$V_g$ 为舱内气体空间体积；$A_t$ 为舱内汽液界面面积；$A_B$ 为破孔处面积；$\rho_l$ 为舱内液体的密度；$\rho_w$ 为水的密度；$u_B$ 为破孔处液体的流速；$C_D$ 为流量系数；$g$ 为引力常数；$a$ 是正弦振荡的振幅；$T$ 为运动周期。

式（6）~（8）相较于式（5）更加全面的分析了油液泄漏入水的过程，当水线下出现泄漏时，将溢油分为两个阶段，第一阶段为由于两侧液体密度不同存在压力差造成的泄漏，这个阶段持续到油舱内容物与周围环境之间通过断裂处及建立流体静水平衡；第二阶段为船舶运动或是波浪造成的压力波动而导致的泄漏，该阶段有水的流入以及水和油液的流出，只要船舶在运动，出流和入流就会一直持续，当破孔两侧被水覆盖后出流和入流都将是水。式（6）~（8）考虑了随着泄漏的不断进行压力发生改变和波浪运动引起的振荡造成的泄漏，更加真实的表述泄漏过程，但是公式中存在很多的系数使得计算过程复杂且很难对系数进行准确确定。根据 C. Soteriou 等人对平口射流研究结论，射流的流量系数主要由泄漏孔内部的流动状态来决定，有三

种流动状态，分别为单相流动、空穴流动以及返流流动，如图 3 所示；根据海洋流体动力学，粗略分析了不同风速下公海波浪的振幅和周期，当风速为 10 节时：$a=0.6m$，$T=2.78s$，当风速 3 节时：$a=0.3m$，$T=1.3s$；在近海区域风速约为 4 节，通常使用的振幅范围为 0.01 m 到 0.05 m，周期范围为 1 s 到 2 s[10-12]。可见，小孔射流模型模拟结果的可信性还需进一步提高，其中存在很多繁琐的系数，且这些系数基本上依靠经验模型得来，导致流量系数、振幅以及周期等参数确定十分困难。

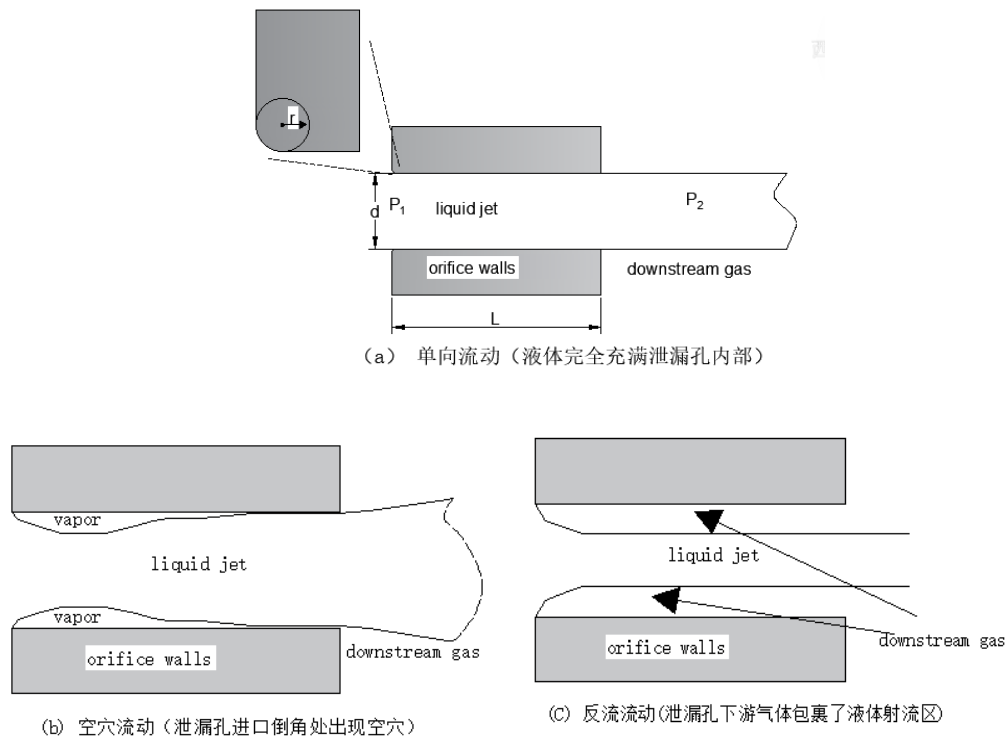

（a）单向流动（液体完全充满泄漏孔内部）

（b）空穴流动（泄漏孔进口倒角处出现空穴）

（C）反流流动(泄漏孔下游气体包裹了液体射流区)

图 3 小孔射流流动状态

上述列出了小孔射流的有关研究成果，这个随液面高度和泄漏时间变化的泄漏源强模型相较于基于溢油量估算的简化模拟更加接近于船舶溢油实际泄漏行为，是目前业界主流的计算模拟方法，计算过程简便占用的计算资源较少，计算的结果更加的精确，但是小孔射流模型模拟计算过程中存在较多的系数，大多数系数依靠经验推到得出，存在一定的不准确性，而且射流模型很难处理复杂油舱形状和不规则破损口面积，对于双层壳船体结构和复杂情景事故模拟也是无能为力，使得小孔射流模型在实际应用中存在很多限制条件，未能被广泛使用[13]。

## 3 基于 CFD 的精确模拟方法

随着计算机的发展，计算流体力学模型（CFD 模型)在理论上能完整地模拟整个溢油过程，是目前认为最能够精确模拟出溢油泄漏源强的方法。计算流体力学模型（CFD 模型)是基于纳维-斯托克斯方程组（N-S 方程）的三维流体守恒定律模型，具有最强的科学依据，能克服射流模型中流量系数难于确定、复杂船体结构无法考虑以及事故风、潮、流情景无法叠加等一系列问题，可以提供油流从油舱通过破损口泄漏入水详细流动行为的描述，也正因为如此，基于 N-S 方程的 CFD 模型被视为研究船舶溢油泄漏行为的最好工具，代表了该领域研究的方向，有不少学者在这方面做了研究。其基本控制方程包括连续性方程、动量方程、能量方程等一系列方程组。

（1）连续性方程（即质量守恒方程）

$$\frac{\partial \rho}{\partial t}+\frac{\partial(\rho u)}{\partial x}+\frac{\partial(\rho v)}{\partial y}+\frac{\partial(\rho w)}{\partial z}=0 \quad (9)$$

其中 $u$、$v$、$w$ 分别为 x、y、z 方向上的速度分量，单位 $m/s$；$\rho$ 为流体密度，单位 $kg/m^3$。

（2）动量守恒方程

$$\frac{\partial(\rho u)}{\partial t}+\frac{\partial(\rho uu)}{\partial x}+\frac{\partial(\rho uv)}{\partial y}+\frac{\partial(\rho uw)}{\partial z}=\frac{\partial}{\partial x}\left(\mu\frac{\partial u}{\partial x}\right)+\frac{\partial}{\partial y}\left(\mu\frac{\partial u}{\partial y}\right)+\frac{\partial}{\partial z}\left(\mu\frac{\partial u}{\partial z}\right)-\frac{\partial P}{\partial x}+S_u \quad (10)$$

$$\frac{\partial(\rho v)}{\partial t}+\frac{\partial(\rho vu)}{\partial x}+\frac{\partial(\rho vv)}{\partial y}+\frac{\partial(\rho vw)}{\partial z}=\frac{\partial}{\partial x}\left(\mu\frac{\partial v}{\partial x}\right)+\frac{\partial}{\partial y}\left(\mu\frac{\partial v}{\partial y}\right)+\frac{\partial}{\partial z}\left(\mu\frac{\partial v}{\partial z}\right)-\frac{\partial P}{\partial y}+S_v \quad (11)$$

$$\frac{\partial(\rho w)}{\partial t}+\frac{\partial(\rho wu)}{\partial x}+\frac{\partial(\rho wv)}{\partial y}+\frac{\partial(\rho ww)}{\partial z}=\frac{\partial}{\partial x}\left(\mu\frac{\partial w}{\partial x}\right)+\frac{\partial}{\partial y}\left(\mu\frac{\partial w}{\partial y}\right)+\frac{\partial}{\partial z}\left(\mu\frac{\partial w}{\partial z}\right)-\frac{\partial P}{\partial z}+S_w \quad (12)$$

式中：$\mu$ 为流体的动力粘度，单位 $N \bullet s/m^2$；$p$ 为作用在流体微元体上的压力，单位 $Pa$；$S_u$、$S_v$、$S_w$ 为动量守恒方程的广义源项。

（3）能量守恒方程

$$\frac{\partial(\rho T)}{\partial t}+\frac{\partial(\rho uT)}{\partial x}+\frac{\partial(\rho vT)}{\partial y}+\frac{\partial(\rho uT)}{\partial z}=\frac{\partial}{\partial x}\left(\frac{k}{c_p}\frac{\partial T}{\partial x}\right)+\frac{\partial}{\partial y}\left(\frac{k}{c_p}\frac{\partial T}{\partial y}\right)+\frac{\partial}{\partial z}\left(\frac{k}{c_p}\frac{\partial T}{\partial z}\right)+S_T \quad (13)$$

式中：$S_T$ 为流体的内热源及由于粘性作用流体机械能转换为热能的部分，$C_p$ 为比热容，$T$ 为温度，$k$ 为流体的传热系数。

通过应用基本控制方程以及研究分析大量的案例，最早提出了相对普适性的预测油轮与另一艘船舶相撞时预期出油量的船舶溢油模型和十二种意外泄漏模型：六个碰撞模型和六个搁浅模型[14-15]，经过长期的实践应用发现存在很大的局限性，又重点研究了船体结构、位移和速度、撞击船位移和速度以及两船的相互作用角对溢油的影响，并提出了新的油轮碰撞和搁浅事故的溢油模型[16]。这些模型的出现简化了溢油事故的分析，能较好的模拟出溢油出流动态，但是已有的都是固定模型，一个模型针对一类溢油事故做出模拟计算，无法做到动态数据的输入，仍然存在很大的局限性。

另一些研究者则以平静的海洋环境为基础建立区域环境，将建好的船体模型导入，再根据实际溢油事故的海况条件输入环境影响因素得到相应的海况条件从而进行模拟研究。早期研究是对破损油舱进行二维条件下建模，如图 4（a）所示，采用 FLUENT 中的 VOF 模型、Segregated、以及压力速度耦合的 PISO 算法来模拟二维条件下的具有自由液面的、不互溶的、油-水-气三相不可压缩非定常流动，能直观的给出溢油的瞬时速度以及溢油的流出动态[17]。后来随着研究进一步的深入，发现二维条件下模型的精确度不能满足实验要求，对模型进行改进，在三维条件下做了溢油的模拟，三维油舱模型如图 4（b）所示，并充分考虑在实际过程中，油品泄漏中油水置换阶段的泄漏行为，深入研究了油的密度、破孔的大小和船舶位于的深度对于溢油速度以及溢油量的影响，并全面研究整个泄漏过程[18]。研究中考虑到船舶在发生意外时的海况是复杂多变的，环境因素对溢油有很大影响，着重分析了在不同风速下的破损处的溢油溢出运动状态，利用 FLUENT 建立数值波浪水槽，模拟波浪运动，建立波浪、水流和风对水下溢油的影响模型，比较分析了不同流速和风速下油液溢出过程[19]。

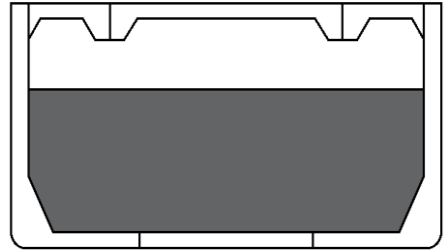

图 4（a）二维油舱模型

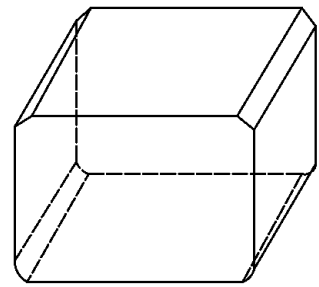

图 4（b）三维油舱模型

上述研究结果表明应用 CFD 比其他方法能更加直观、准确的模拟出每时刻溢油的动态，能更好的结合环境因素对溢油的影响作用，是该领域的研究发展方向。CFD 模型可以近似模拟出溢油动态过程，在理想情况下理论上能够模拟出溢油的全部过程，但是由于油舱模型与实际油舱结构存在一定差别以及实际中海况条件的变化、多参数如何叠加、计算方法选定上还在进一步研究当中，需要不断改进和完善。

# 4 模拟方法改进方向

近几十年来，世界上的许多国家机构和研究者在船舶溢油泄漏模拟上做了大量深入的研究，相继提出了各自的理论成果，取得了一定的成就。当然，现阶段研究成果还存在许多不完善的地方，但是它已经为溢油的应急处理和危害鉴定评估提供了很大的支持作用，相信随着研究的不断深入，更加可靠、精确的溢油模型将会建立。基于计算流体力学的数值模拟是近年来研究发展的方向，具有最强的科学依据，特别是船舶溢油油液泄漏入水行为过程这一基础数据的模拟研究，这将是后续准确研究的前提，但是，目前的 CFD 模拟未能考虑到复杂船体结构对溢油的影响和综合分析风、流、浪、航速、航向等事故情景参数与多相流之间的影响关系以及离散方式和求解方法的设定。因此，对未来进一步研究工作的建议是：

溢油模型的构建——复杂结构下的船舶模型。目前的研究是对单个油舱进行建模，分析油舱在各种影响因素作用下的油液泄漏情况，而实际中船体外壳结构复杂以及油舱大小和布置位置各有不同，这些都会对溢油的流速、流动状态等产生影响，因此需要对整个船体进行建模，这样才能更加真实的进行模拟研究。今后研究中应重点考虑各种舱室配置设计，包括舷侧纵舱壁、中纵舱壁、横舱壁等，如图 5 所示。复杂油舱结构虽然不能够完全消除油舱损坏造成的污染风险，但复杂油舱结构对于油液溢出有很大影响作用。双壳结构在发生轻微碰撞时可以防止内壳破裂，海水会进入压载舱不会造成油液泄漏；发生严重碰撞导致内壳破损时，双层结构会截留一部分油液以及对油液泄漏速度和流动状态产生一定影响。因此，深入研究油舱复杂结构对于油液泄漏入水过程模拟的准确性研究有重要意义。

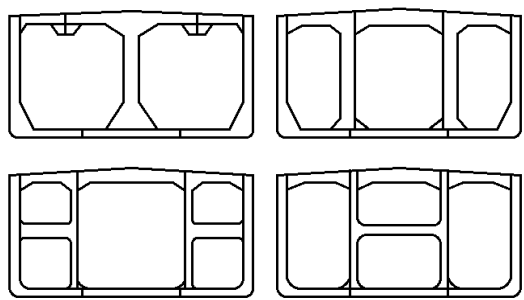

图 5 油舱布置简图

事故情景参数的设定，影响溢油速度和溢油量的因素有很多，各种影响因素是相互制约、相互影响的。目前的数值模拟是基于单因素作用下的溢油模拟分析，但是溢油事故的发生通常是在恶劣环境条件下，且多种溢油影响因素同时存在，因此利用 FLUENT 软件建立三维数值波浪水槽，模拟波浪运动，建立波浪、水流和风作用下的水下溢油模型，综合分析风、流、浪、航速、航向等事故情景参数对溢油油液泄漏入水动态变化的影响。

离散方式和求解方法设定，深入研究 CFD 仿真基本原理和计算软件，比较分析 Eulerian-Lagrangian（欧拉－拉格朗日）和 Eulerian-Eulerian（欧拉－欧拉）等多相流计算方法优缺点，结合船舶溢油泄漏行为特征初步分析，确定 CFD 多相流数值计算方法适用条件，在不同阶段采用不同的计算方法，充分发挥每个求解方法的优点，使得计算结果更加准确。